\documentclass[12 pt]{amsart}
\begin{document}

\title{Toric Varieties in Hilbert Schemes}
\author{Heather Russell}
\maketitle

\tableofcontents
\date{\today}
\def\Q{{\mathbb Q}}
\def\A{{\mathbb A}}
\def\P{{\mathbb P}}
\def\C{{\mathbb C}}
\def\G{{\mathbb G}}
\def\Z{{\mathbb Z}}
\def\Sp{{\textup{Spec}}}
\def\O{{\mathcal O}}
\def\H{{\textup{Hilb}}}
\def\N{{\mathbb N}} 
\def\I{{\mathcal I}}
\def\m{{\mathfrak m}}
\def\Se{{\mathcal S}}
\def\aut{{\rm Aut}}

\newtheorem{example}{Example}[section]
\newtheorem{theorem}{Theorem}[section]   
\newtheorem{corollary}{Corollary}[section]
\newtheorem{question}{Question}[section]
\newtheorem{lemma}{Lemma}[section]
\newtheorem{proposition}{Proposition}[section]

\section{Introduction} Given a field $K$, let $R$ be the 
power series ring $K[[x_1, \dots, x_n]]$.  There is a natural 
action of $\aut (R)$ on the Hilbert scheme 
$\rm{Hilb}^d(R)$ parameterizing 
ideals in $R$ of colength $d$.  Recall from [10] that given a
smooth variety $X$ of dimension $n$ and a monomial ideal $I \subset R$
of finite colength, the space $U(I)$ parametrizing subschemes of $X$ 
isomorphic to $\Sp (R/I)$ sits naturally inside a suitable Hilbert
scheme $H$.  In [10] a flag bundle $B$ on the tangent
bundle of $X$ and a space $Y \subset B \times U(I)$ such that 
$Y$ is a fiber bundle over $B$ with respect to the first projection
map and an \'{e}tale covering of $U(I)$ with respect to the second is 
constructed.  Moreover, $Y$ 
has the property that its closure in $B\times H$ is also a fiber bundle 
over $B$.   Via an isomorphism $\varphi$ of $R$ and the
local ring of the point of $X$ over which the fiber lies, 
the fiber is isomorphic to the  
closure of the $G$-orbit of $I$ in $H$ where $G$ is
the subgroup of $\aut (R) $ fixing the flag over which the fiber lies
via 
$\varphi$.  The complexity of the fiber $F(I)$ of $Y$ over $B$ and its
closure is partially 
measured by the measuring sequence of $I$ as defined in [10].  
In the first level of complexity, the $G$ equivariance of $\bar F (I)$
forces its normalization to be a projective space.  In the second
level, the normalization is forced to be a toric variety.  In the
third level, we lump together all of the other cases.  The following 
three sections correspond to these three levels of difficulty, each 
illustrating techniques for understanding $F$ and $\bar F (I)$.  For
simplicity, we restrict ourselves to the case when $n =2$, thus
freeing the variable $n$ for later use.  Moreover, we will use the
variables $x$ and $y$ in place of $x_1$ and $x_2$, to cut down on the
use of subscripts.  In this case it is very easy to define the
measuring sequences.  

\vspace{.1 in}

\noindent{\bf Definition:} Say a sequence of monomial ideals 
$I_1, \dots I_r$ of $R$, each having finite colength, has 
{\it measuring sequence} $(x,y^a), (x^b,y)$ denoted $m(a,b)$ 
where $a$ and $b$ are the smallest
positive integers such that any autormorphism of $R$ fixing the
measuring sequence also fixes the original sequence of ideals.  We
will say that $m(a,b) \le m(c,d)$ if $a\le c$ and $b \le d$.  

\vspace{.1 in}

If the measuring sequence of $I$ is $m(1,1)$ then $I$ must be a power
of the maximal ideal $\m$ of $R$.  Otherwise, up to permuting $x$ and
$y$, the group $G$ is the group of automorphisms sending the ideal 
$(x,y^2)$ to itself and the fiber $F(I)$ of $Y$ over $B$ is isomorphic
to the quotient of $G$ by the subgroup $G(I)$ of
automorphisms fixing $I$.   
To avoid defining $G$ repeatedly, we will henceforth let $G$ be the
group of automorphisms fixing $(x,y^2)$.  The measuring sequence of
the ideal 
$I$ determines the group $G(I)$.  More generally, if we consider the $G$
orbit of a sequence of monomial ideals $F(I_1, \dots, I_r)$, then this
space is isomorphic to the quotient of $G$ by its subgroup $G(I_1, \dots,
I_r)$ of automorphisms fixing the $I_j$'s.  The group 
$G(I_1, \dots, I_r)$ is likewise determined by the measuring sequence
of the $I_j$'s.     
     
Given a monomial ideal $I \in R$ let $\Gamma (I)$ denote the
staircase associated to $I$.  Then $\Gamma (I)$ consists of the unit 
squares tiling the plane with lower left hand corner giving the exponent
vector for a monomial in the complement of $I$.  In characteristic
zero, the measuring sequence of ideals $I_1, \dots ,I_r$
is of the form $m(a,b)$ where $a$ is the largest
number of boxes by which two consecutive rows differ in any of the 
$\Gamma (I_j)$'s and and $b$ the
largest number by which two consecutive columns differ.

\section{Finding the Boundary as a Limit Point}

If the sequence of ideals $I$ has measuring sequence at
most $m(2,1)$, then $F(I)$ is a single point.  If it is 
$m(3,1)$ or $m(2,2)$, then the space 
$F(I)$ is an affine line and is closure 
$\bar F(I)$ is a $\P^{1}$.  Thus the boundary is a single
point.  In this section we give an algorithm for finding the boundary
point simply by shifting boxes of the staircases $\Gamma_i$.  

By Lemma 3.3 of [10], 
if the measuring sequence of $I$ is $m(3,1)$ (respectively $m(2,2)$) 
a set of coset representatives for 
$G/G(I_1, \dots ,I_r)$ are given by automorphism of the form $g(t)$
where $g(t)(x)= x+ty^2$ and $g(t)(y)= y$ (respectively $g(t)(x) = x$ 
and $g(t)(y) = y+tx$).  Thus the boundary point is the limit as $t$
goes to infinity of $g(t)(I)$ and can be found by 
Proposition~\ref{limits} below.

\vspace{.1 in} 
                  
\noindent {\bf Definition:}  If the characteristic $p$ of $K$ is positive, 
let $\prec$ be the
total ordering on the set of non-negative 
integers defined by  $a \prec b$ if either 
${\rm ord} _p (a) > {\rm ord} _p (b)$ or 
${\rm ord} _p (a) = {\rm ord} _p (b)$and   
$a < b $.  Otherwise, let it be the usual ordering.  
Given a finite collection of non-negative 
integers $T$, let $f$ be the injective map from $T$ to the set of non-negative 
integers sending each element of $T$ in ascending order with respect to 
$\prec$ to the smallest integer divisible by a greater or equal power of $p$.  
Define the $p$-shift $T^{\prime}$ of
$T$ to be the image of $f$.   
 
\begin{lemma}\label{limlem} Let $T$ be a finite collection of non-negative 
integers.  The limit $V(T)$ as $t$ goes to infinity   
of the subspaces $V(T,t)$ of 
$K[x]$ generated by the polynomials $(x+t)^e$ for exponents $e
\in T$ is spanned by the monomials $x^e$ for exponents $e \in T^{\prime}$
where $T^{\prime}$ is the $p$-shift of $T$.  
\end{lemma} 

\noindent {\bf Proof:} 
Let $r = |T|$ and assume that the lemma holds for all 
smaller cardinalities.  Suppose that $T^{\prime}$ is the 
set of integers from $0$ to $r-1$.  Let $M$ be the span of wedge products 
of monomials $\wedge ^r(K[x])$ such that for some power of $p$ there are more 
exponents of the monomials divisible by that power than elements of 
$T^{\prime}$ divisible by that power.  Then $M$ is fixed by $g(t)$ for 
any $t$.  Thus so is the limit as $t$ goes to infinity of the images of 
the projections of the spaces $V(T,t)$ in $\P (\wedge ^r (K[x])/M)$.   
The projection map is well-defined because the coefficient 
of the wedge product of monomials with exponent in $T$ is non-zero in the 
image of $V(T,t)$ in $\P (\wedge ^r(K[x])$.  Since the only 
$g(t)$-fixed point in $\P (\wedge ^r (K[x])/M)$ 
is the wedge product of monomials with exponent in $T^{\prime}$, 
it follows that the coefficient of the wedge product of 
monomials with exponent in $T^{\prime}$ in $\P (\wedge ^r(K[x])$ is 
non-zero and has the highest degree in $t$.  Thus the limit $V(T)$ is 
as claimed.  
         
Suppose that $T^{\prime}$ is not the 
set of integers from $0$ to $r-1$. Let $f:T \rightarrow T^{\prime}$ be the 
map taking a label of a box to the label of the box it shifts to as in the 
definition of $p$-shift above.  
Let $T_1$ and $T_2$ be the compliments in the preimage under $f$ of the 
compliments of the largest element of $T^{\prime}$ and the largest element 
of $T^{\prime}$ with the smallest power of $p$ dividing it.  Then 
$$T_1^{\prime} + T_2^{\prime} = T^{\prime}$$     
and thus 
$$V(T_1) + V(T_2) \subset V(T).$$
It follows from the inductive hypothesis applied to $T_1$ and $T_2$ 
that $V(T)$ is as claimed.    

\vspace {.1 in}

\begin{proposition}\label{limits} Let $I\subset R$ be 
a monomial ideal of finite colength.  Let $h$ be a monomial not divisible by 
$x_i$.  For each $t \in K$, let $g(t)$ be the element of 
${\aut}(R)$ with 
$$g(t)(x_j)= x_j$$ for $j \ne i$ and 
$$g(t)(x_i)= x_i + th.$$  
Let $J$ be the flat limit of ideals 
$$\lim _{t \to \infty}g(t)(I).$$ 
Then $J$ is characterized by the following.  
Let $T$ be the set of integers $a$ such that 
a basis of $M$ is given by monomials of the form 
$x_i^ah^{d-a}f$ for $a \in T$ and $f$ a monomial of minimal
degree.  Then the corresponding graded piece of $J$ has a basis of
monomials of the form $x_i^ah^{d-a}f$ for $a \in
T^{\prime}$.   
\end{proposition}

\noindent{\bf Proof:}  Since $f$ and $h$ are invariant under 
$g(t)$, the proof reduces to Lemma~\ref{limlem}.    

\vspace{.1 in}

\section{Toric Varieties as Fibers}      

If the measuring sequence is 
$m(4,1)$ or $m(2,3)$ by Lemma 3.3 [10], coset representatives for
$G/G(I_1, \dots ,I_r)$ are given by the set of automorphisms 
$g(a,b)$ with $g(a,b) = x + ay^2 + by^3$ and $g(a,b)(y)= y$ and the
set of automorphisms $h(a,b)$ with $h(a,b)(x)= x+ay^2$ and $h(a,b)(y)=
y+bx$ respectively where $a$ and $b$ range over $K$.  Hence 
$F(I_1, \dots ,I_r)$ is an affine plane with $a$ and $b$ as
coordinates.  Its closure, being equivariant under the action of $G$
is equivariant under the action of 
scaling $x$ and $y$.   
One can find the fan of the
normalization of $\bar F(I_1, \dots ,I_r)$ as follows.  Let $V_j$ be
the quotient of two monomial ideals between which all ideal in
$F(I_j)$ are sandwiched.  Then $F(I_j)$ has an embedding in the
Grassmanian of subspaces of $V_j$ of the appropriate dimension which
can in turn be embedded in projective space by Pl{\"u}cker coordinates.
Taking the coordinates to correspond to wedge products of monomials, 
the coordinate functions will be monomials in $a$ and $b$ [10].  
Since the product of these projective spaces for each $j$ can be
embedded into one big projective space by the Segre embedding, 
coordinate functions embedding $F(I_1, \dots ,I_r)$ in this big
projective space are given by products of the coordinate functions 
embedding each of the $F(I_j)$'s, one for each $j$.  Since these will
also be monomials, we can plot the exponent vectors $(m,n)$ of each
coordinate function $a^mb^n$ and take the
convex hull $H$.  The cone covering the plane with a ray normal to
each edge of $H$ is a fan $\Delta (I_1, \dots ,I_r)$ 
for the normalization $F(I_1, \dots ,I_r)$ in the sense of \cite{F}.  
  
\vspace{.1 in}

\noindent {\bf Definition:} Given a sequence of ideals 
$(I_1, \dots,I_r)$ with measuring sequence $m(4,1)$ or $m(3,2)$, 
say that the fan $\Delta (I_1, \dots ,I_r)$ constructed above is the 
{\it standard fan} for $F(I_1, \dots ,I_r)$.  

\vspace{.1 in}

Suppose one wanted to find the standard fan for the space
$F(I(4)^m)$.  The ideals in the $G$ orbit of $I(4)^m$ are sandwiched
between $I(2)^m$ and $I(2)^{2m}$.  Let $x$ have weight $2$ and $y$
have weight $1$.  Then the $2m^2$ polynomials of the form $g(a,b)(x^c)x^dy^e$
for $4c + 2d+ e$ equal to either to $4m$ or 
$4m+1$ and $c>0$ span $I(4)^m$ as a subspace of $I(2)^m/I(2)^{2m}$.  
Let $M_m$ be the matrix with rows corresponding to these
generators arranged in increasing order first by $c$ and then by $d$,
columns corresponding to the monomials spanning $I(2)^m/I(2)^{2m}$
arranged in order of descending weight and descending powers of $y$,
and entries given the coefficient of the column's monomial in the row's 
polynomial.  Coordinate functions embedding $F(I(4)^m)$ into
projective space are given by the determinants of the maximum minors of 
$M_m$.  With no additional insight this becomes a hefty task as $m$
becomes large.  However, with the help of a few observations we will 
find the standard tori for $F(I(4)^m)$ without taking a single
determinant.  Letting  $x$ have weight $(1,0)$, $y$ have weight $(0,1)$
respectively, $a$ have weight $(1,-2)$ and $y$ have weight 
$(1,-3)$.  Then $g(a,b)(x)$ is homogenious and thus the wedge product
of the $2m^2$ polynomials are homogenious and so the sum of the
weights of the determinant of a maximum minor and the monomials
corresponding to the columns is constant.  Thus, the powers of $a$ and
$b$ in the determinant of a minor can be read off easily, but it can
be difficult to tell even whether the coefficient is zero or not.  To
find the standard fan of $I(4)^m$, it is only necessary to find the
vertices of the convex hull of the exponent vectors of the non-zero
coordinate functions.  These vertices correspond to the monomial
ideals in $F(I(4)^m)$ in the sense that the one non-zero coordinate 
of the monomial ideal maps to the corresponding vertex.  Thus we need
only check the coefficients of coordinates corresponding to monomial
ideals, but this may still be unneccesarily difficult.  For
example, the coefficient of the 
coordinate function corresponding to the ideal 
$$I(\underbrace{1,\dots ,1}_{m},2,\underbrace {1,\dots, 1}_{m-1})$$ 
is the determinant of the $2m^2 \times 2m^2$ 
matrix which can be described as follows.  Divide up the rows an
columns as follows.  Put the first $2m-1$ rows together, then the next
$2n-1$ rows together, then the next $2m-3$ rows together, then the
next $2n-3$ rows together and so on, alternately keeping the number of
rows the same and decreasing the number of rows by two.  Divide the
columns up similarly except that the first group of columns will have
one extra column and the last group will have one less and hence no
columns.  Fill in the rectangles as follows.  In the $(2n-1)^{\rm st}$
and $(2n)^{\rm th}$ 
row of rectangles, the entry in the $i^{\rm th}$ and $j^{\rm th}$ column
of the $k^{\rm th}$ rectangle is $$\left ( \frac{n}{k+1-n} \right )
\left ( \frac {k+1-n}{j-i} \right ) ~~ {\rm and } ~~ 
\left ( \frac{n}{k-n} \right )
\left ( \frac {k-n}{j-i} \right )$$ 
respectively.  
For example, the matrix $M_3$ is shown in Figure~\ref{matrix}. 

\begin{figure}\label{matrix}
$$\left [   
\begin{array}{cccccc|ccccc|cccc|cc|c}
1&1&0&0&0&0& 0&0&0&0&0& 0&0&0&0& 0&0 &0\\
 0&1&1&0&0&0& 0&0&0&0&0& 0&0&0&0& 0&0 &0\\
 0&0&1&1&0&0& 0&0&0&0&0& 0&0&0&0& 0&0 &0\\
 0&0&0&1&1&0& 0&0&0&0&0& 0&0&0&0& 0&0 &0\\
 0&0&0&0&1&1& 0&0&0&0&0& 0&0&0&0& 0&0 &0\\ \hline 
 1&0&0&0&0&0& 1&1&0&0&0& 0&0&0&0& 0&0 &0\\
 0&1&0&0&0&0& 0&1&1&0&0& 0&0&0&0& 0&0 &0\\
 0&0&1&0&0&0& 0&0&1&1&0& 0&0&0&0& 0&0 &0\\
 0&0&0&1&0&0& 0&0&0&1&1& 0&0&0&0& 0&0 &0\\
 0&0&0&0&1&0& 0&0&0&0&1& 0&0&0&0& 0&0 &0\\ \hline
 1&0&0&0&0&0& 2&2&0&0&0& 1&2&1&0& 0&0 &0\\
 0&1&0&0&0&0& 0&2&2&0&0& 0&1&2&1& 0&0 &0\\
 0&0&1&0&0&0& 0&0&2&2&0& 0&0&1&2& 0&0 &0\\ \hline 
 0&0&0&0&0&0& 1&0&0&0&0& 2&2&0&0& 1&2 &0\\
 0&0&0&0&0&0& 0&1&0&0&0& 0&2&2&0& 0&1 &0\\
 0&0&0&0&0&0& 0&0&1&0&0& 0&0&2&2& 0&0 &0\\ \hline 
 0&0&0&0&0&0& 1&0&0&0&0& 3&3&0&0& 3&6 &1\\ \hline 
 0&0&0&0&0&0& 0&0&0&0&0& 1&0&0&0& 3&3 &3
\end{array}
 \right ]$$
\caption{ $M_3$}
\end{figure}

We will adopt the following notation.  Given an ideal $I$ with 
measuring sequence $m(4,1)$, let $I^+(m,n)$ (respectively
$I^-(m,n) $) be the ideal corresponding to the cone just clockwise
(respectively counterclockwise) to the ray through $(m,n)$ in the
standard fan $\Delta (I)$.  We extend
this definition to ideals with smaller measuring sequence in the
natural way.  Similarly, if $I$ has measuring sequence 
$m(3,2)$ we let $I_+(m,n)$ (respectively $I_-(m,n)$) 
be the ideal corresponding to the cone just clockwise
(respectively counterclockwise) to the ray through $(m,n)$ in the
standard fan $\Delta (I)$.
and we extend this definition to ideals with smaller measuring sequence in
the natural way.

\begin{proposition}\label{endpt} 
Given an ideal $I$ with 
measuring sequence $m(4,1)$ or $m(3,2)$, 
the standard fan $\Delta (I)$ has lower left quadrant 
corresponding to the affine plane $F (I)$ with the negative $x$ and $y$ 
axes corresponding to the line $a=0$ and $b=0$ respectively.  In the
first case 
$$I^+(-1,0) = 
I^-(0,1)= 
\lim _{t \to \infty}g(0,t)(I),$$
$$I^-(0,-1) =  
\lim _{t \to \infty}g(t,0)(I),$$
$$I^+(0,1) = 
\lim _{t \to \infty}g(t,0)I^-(0,1),$$
and, if the characteristic of $K$ is not $2$, then 
$$I^-(0,-1) =
I^+(1,2)$$
and $$I^-(1,2) = 
\lim _{t \to \infty}g(0,t)I^+(1,2).$$
In the second case, 
$$I_+(-1,0)= 
I_-(0,1) = 
\lim _{t \to \infty}h(0,t)(I),$$
$$I_-(0,-1) = 
I_+(1,0)=
\lim _{t \to \infty}h(t,0)(I),$$
$$I_+(0,1) = 
\lim _{t \to \infty}h(t,0)I_-(0,1),$$
and 
$$I_-(0,1) = 
\lim _{t \to \infty}h(0,t)I_+(0,1).$$

\end{proposition}

\noindent{\bf Proof:}  One can verify by hand that 
$$(x,y^3)^+(-1,0)= (x,y^3)^-(0,1) = (x,y^3),$$
$$(x,y^4)^-(0,-1)= (x,y^4)^+(1,2)= (x^2,y^2),$$
$$(x,y^3)_+(-1,0)= (x,y^3)_-(0,1) = (x,y^3),$$
$$(x^2,y)_-(0,-1)= (x^2,y)^+(1,0)= (x^2,y).$$
Except for the ideal in the second line in characteristic $2$, 
these ideals are not invariant under the action of $G$.  Thus 
the proposition follows.    

\vspace{.1 in}

See \cite{Col} and [10] for the relevant definitions in the
following corollary.  

\begin{corollary} The Semple bundle $F(4)$ is isomorphic to 
the alignment correspondence $C((x,y^2), (x,y^3), (x^2,xy,y^5))$.  
\end{corollary}

\noindent{\bf Proof :} Recall that the Semple bundle $F(4)$ is an 
compactification of the space of curvilinear $3$-jets on $X$ with the 
property that it is a fiber bundle over the projectivized 
tangent bundle of $X$ with fiber equivariant under the action of $G$.
Moreover, the fibers are smooth and have two boundary divisors with 
self-intersection $0$ and $-3$ respectively, as one can check from 
the Chow ring of $F(4)$ computed in \cite{Col}.  Thus, 
they are toric varieties 
and giving the interior of the fiber in a manner compatible with 
the standard fans constructed before, the only possible fan for the 
fiber of $F(4)$ is the fan for $F((x,y^3), (x^2,xy,y^5))$.  

\begin{proposition}\label{slopes} Given an ideal $I$ 
with measuring sequence $m(4,1)$ (respectively $m(3,2)$), let $I[n,m]$ be
the limit as $t$ tends to infinity of $g(at^n,bt^m)(I)$ (respectively 
$h(at^n,bt^m)$).  As $\frac{a^n}{b^m}$ varies 
through $\C ^*$, $I[n,m]$ runs over the torus orbit corresponding to
the ray through $(m,n)$ in the standard torus $\Delta (I)$.  
The limits as $a$ and $b$ tend to
$0$ of $I[m,n]$ are $I^+(m,n)$ and $I^-(m,n)$ (respectively  
$I_+(m,n)$ and $I_-(m,n)$).  In particular, 
giving $x$ weight $n-m$ (respectively n+m) and $y$ weight $2n-3m$ 
(respectively $n+2m$),  
each graded piece of the monomial ideal 
$I^+(m,n)$ (respectively $I_+(m,n)$) has the same dimension as the respective
graded piece of the monomial ideal $I^-(m,n)$ respectively ($I_-(m,n)$). 
\end{proposition}

\noindent{\bf Proof:} By the construction of the standard torus 
$\Delta (I)$, if there is a ray in $\Delta (I)$ through $(m,n)$, 
the limit of a path $(at^n,bt^m)$ with $ab \ne 0$ in the affine
coordinates of $F(I)$, as $t$ goes to infinity, lands in the torus
orbit corresponding to that ray, the points of the torus orbit, being
in bijection with the ratios $\frac{a^n}{b^m}$.  Thus as $a$ and $b$
tend to $0$, we get the two ideals 
$I^-(m,n)$ and $I^+(m,n)$ (respectively $I_+(m,n)$ and $I_-(m,n)$).  
The point
$(at^n,bt^m)$ corresponds to the ideal $g(at^n,bt^m)(I)$ 
(respectively $h(at^n,bt^m)(I)$).  Giving $a$ weight $-n$ 
and $b$ weight $-m$, the elements of $I[m,n]$ are the coefficients
of the highest powers of $t$ in elements of $g(at^n,bt^m)(I)$ 
(respectively $h(at^n,bt^m)(I)$).  Thus $I[m,n]$ is generated by homogenious 
elements.  Moreover giving 
$x$ weight $n-m$ (respectively n+m), $y$ weight $2n-3m$ 
(respectively $n+2m$), and keeping the weights of $a$ and $b$ the
same, $g(at^n,bt^m)(I)$ (respectively $h(at^n,bt^m)(I)$) 
is generated by homogenious elements.  Thus  $I[m,n]$, being already
homogenious in $a$ and $b$ is homogenious in $x$ and $y$ with respect
to the weights.  Therefore the dimensions of the graded pieces of the 
two limits $I^-(m,n)$ and $I^+(m,n)$ (respectively $I_+(m,n)$ and
$I_-(m,n)$) with respect to the weights on
$x$ and $y$ are the same as the dimensions of
the respective graded pieces of $I[m,n]$.

\begin{corollary}\label{slopecor} Given a monomial ideal $I$ 
with measuring sequence $m(4,1)$ (respectively $m(3,2)$) 
if two monomial ideals $I_1$ and $I_2$ 
correspond to boundary points of $F(I)$ and the quotient of the
product of monomials in there respective complements is
$\frac{x^c}{y^d}$ with $c > 0$, then the ray $r$ through
the point $(2c-d,3c-d)$ (respectively $(2c-d, d-c)$) lies in the 
convex cone bounded by the two cones in the standard fan 
$\Delta (I)$ corresponding to 
$I_1$ and $I_2$.  
\end{corollary}

\noindent{\bf Proof:}  If $I_1$ and $I_2$ correspond to adjacent cones
in $\Delta(I)$, then by Proposition~\ref{slopes}, $r$ is the ray
between them.  The rest follows from elementary properties of
medians.  

\begin{proposition}\label{mult}  Given ideals $I_1$, $I_2$, and $I_3$ with 
$I_1I_2 \subset I_3$ then for any point $P \in F(I_1, I_2, I_3)$, 
we have $\pi _1(P)\pi_2(P)\subset \pi_3(P)$ where $\pi_i$ is 
projection to $F(I_i)$.
\end{proposition}

\noindent{\bf Proof:}  Let $V_1$, $V_2$, and $V_3$ be quotients of 
ideals between which the ideals in the $G$ orbits of $I_1$, $I_2$, and
$I_3$ lie respectively.  Consider the incidence correspondence
$$\Gamma = \{ (\alpha_1, \alpha_2, \alpha_3, (a_1, a_2, a_3))\in 
V_1\times V_2 \times V_3 \times \bar F(I_1, I_2, I_3): $$
$$\alpha_i \in a_i, \alpha_1\alpha_2 = \alpha_3 \}.$$  Then the set 
$$\{ (\alpha_1, \alpha_2, \alpha_3, (a_1, a_2, a_3))\in 
V_1\times V_2 \times V_3 \times F(I_1, I_2, I_3): \alpha_i \in a_i\}$$ 
is contained in $\Gamma$ and hence since so is its closure  
$$\{ (\alpha_1, \alpha_2, \alpha_3, (a_1, a_2, a_3))\in 
V_1\times V_2 \times V_3 \times \bar F(I_1, I_2, I_3): \alpha_i \in a_i\}.$$

\vspace{.1 in}

Given a sequence of non-negative 
integers,
$$s= n_1, \dots, n_r,$$ we will let $I(s)$ denote the ideal 
$$I(s)= (x^r, x^{r-1}y^{n_1}, \dots , y^{n_1+\dots +n_r}).$$

We will use the following two observations repeatedly without
mentioning so explicitely. 

\begin{enumerate} 
\item $${\rm col}(I(n_1, \dots, n_r))= \sum_{i=1}^r (r+1-i)n_i.$$
\item $I(n_1, \dots, n_r)I(m_1, \dots, m_k)= I(l_1, \dots, l_{r+k})$
  where $$l_i = \min_{\alpha +\beta = i}\left( \sum_{j=1}^{\alpha} n_j
  = \sum_{j=1}^{\beta} m_j \right ).$$
\end{enumerate}

We will use the notation $$I^{\pm}(r,(m,n))=
((x,y^4)^r)^{\pm}(m,n).$$  From the results proved in this section, we
will now demonstrate the following eight claims.

\begin{figure}\label{boundary}
\begin{picture}(350,450)
\put (270,0){\line(0,1){40}}
\put (270,20){\small $(0,1)$}
\put (230,40){\tiny $2,2,2,2,1,2,2,2,1$}
\put (270,50){\line(0,1){40}}
\put (270,70){\small $(1,4)$}
\put (230,90){\tiny $1,2,2,1,2,1,2,1,2,1$}
\put (270,100){\line(0,1){40}}
\put (270,120){\small $(1,3)$}
\put (230,140){\tiny $2,1,2,1,2,1,2,1,1,2$}
\put (270,150){\line(0,1){40}}
\put (270,170){\small $(3,8)$}
\put (230,190){\tiny $1,1,2,1,1,2,1,1,2,1,1$}
\put (270,200){\line(0,1){40}}
\put (270,220){\small $(5,12)$}
\put (230,240){\tiny $1,2,1,1,1,2,1,1,1,2,1$}
\put (270,250){\line(0,1){40}}
\put (270,270){\small $(7,16)$}
\put (230,290){\tiny $2,1,1,1,1,2,1,1,1,1,2$}
\put (270,300){\line(0,1){40}}
\put (270,320){\small $(9,20)$}
\put (230,340){\tiny $1,1,1,1,1,1,2,1,1,1,1,1$}
\put (270,350){\line(0,1){40}} 
\put (270,370){\small $(1,2)$}
\put (265,400){6.}
  
\put (190,30) {\line(0,1){40}}
\put (190,50) {\small $(0,1)$}
\put (170,70) {\tiny $1,2,2,2,1,2,2,2$}
\put (190,80) {\line(0,1){40}}
\put (190,100){\small $(1,4)$}
\put (170,120){\tiny $2,2,1,2,1,2,1,2$}
\put (190,130){\line(0,1){40}}
\put (190,150){\small $(1,3)$}
\put (170,170){\tiny $1,1,2,1,2,1,2,1,1$}
\put (190,180){\line(0,1){40}}
\put (190,200){\small $(3,8)$}
\put (170,220){\tiny $1,2,1,1,2,1,1,2,1$}
\put (190,230){\line (0,1){40}}
\put (190,250){\small $(5,12)$}
\put (170,270){\tiny $2,1,1,1,2,1,1,1,2$}
\put (190,280){\line (0,1){40}}
\put (190,300){\small $(7,16)$}
\put (170,320){\tiny $1,1,1,1,1,2,1,1,1,1$}
\put (190,330){\line (0,1){40}}
\put (190,350){\small $(1,2)$}
\put (185,380){5.}

\put (135,60) {\line (0,1){40}}
\put (135,80) {\small $(0,1)$}
\put (115,100) {\tiny $2,2,2,2,1,2$}
\put (135,110){\line (0,1){40}}
\put (135,130){\small $(1,4)$}
\put (115,150){\tiny $1,2,1,2,1,2,1$}
\put (135,160){\line (0,1){40}}
\put (135,180){\small $(3,8)$}
\put (115,200){\tiny $2,1,1,2,1,1,2$}
\put (135,210){\line (0,1){40}}
\put (135,230){\small $(5,12)$}
\put (115,250){\tiny $1,1,1,1,2,1,1,1$}
\put (135,260){\line (0,1){40}}
\put (135,280){\small $(1,2)$}
\put (130,310){4.}

\put (87,90) {\line (0,1){40}}
\put (87,110){\small $(0,1)$}
\put (67,130){\tiny $1,2,2,2,1$}
\put (87,140){\line (0,1){40}}
\put (87,160){\small $(1,4)$}
\put (67,180){\tiny $2,1,2,1,2$}
\put (87,190){\line (0,1){40}}
\put (87,210){\small $(3,8)$}
\put (67,230){\tiny $1,1,1,2,1,1$}
\put (87,240){\line (0,1){40}}
\put (87,260){\small $(1,2)$}
\put (82,290){3.} 

\put (25,160){\tiny $2,2,2$}
\put (40,170){\line (0,1){40}}
\put (40,190){\small $(1,4)$}
\put (25,210){\tiny $1,1,2,1$}
\put (40,220){\line (0,1){40}}
\put (40,240){\small $(1,2)$}
\put (35,270){2.}

\put (0,190) {\tiny $1,2$}
\put (5,200){\line (0,1){40}}
\put (5,220){\small $(1,2)$}
\put (0,250){1.}

\end{picture}
\caption{Boundary Components}
\end{figure}
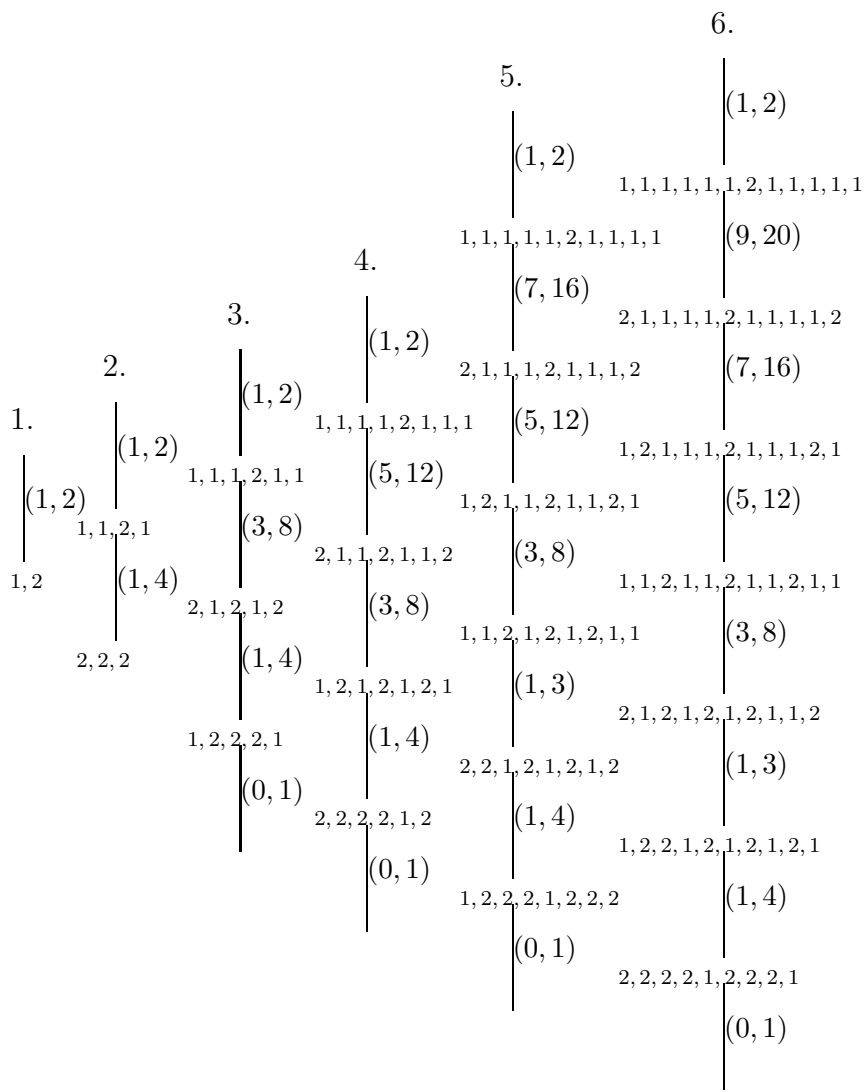

\begin{enumerate}
\item $I^-(n,(1,0)= I^+(n,(1,2))= (x^2,y^2)^n.$

\item $I^-(n, (1,2))= 
I(\underbrace{1, \dots ,1}_n,2,\underbrace{1, \dots ,1}_{n-1})= 
I^+(n, (2n-3,4n-4))$

\item Letting $s_1 = 1,2$, 
$s_2 = 2,2,2$ and $s_3 = 3,2,2,2$ then 
$$I^+(3q+r,(-1,0))=I(s_r,\underbrace{s_3, \dots ,s_3}_q)=
I(s_r)I(s_3)^q$$ for $1\le r \le 3.$ 

\item Letting $s_1 = 2,2,1$, $s_2 = 2,2,2$ and $s_3 = 2,1,2$,   
$$I^+(2n+1,(0,1))= I^-(2n+1,(1,4))=
I(1,2, \underbrace {s_1, s_2, s_3, \dots }_n)$$ and 
$$I^+(2n,(0,1))= I^-(2n,(1,4))
=I(\underbrace{s_2, s_3, s_1, \dots }_n).$$   

\item If  $\left [ \frac {n}{2} \right ] \le k \le n-1$ then 
$$I^-(n, (2k-1,4k))= 
I(\underbrace{1, \dots ,1}_{n-k-1},2,\underbrace{1, \dots ,1}_{k-1},2,
\underbrace{1, \dots ,1}_{k-1},2,\underbrace{1, \dots ,1}_{n-k-1}).$$

\item  If  $\left [ \frac {n-1}{2} \right ] \le k \le n-2$ then 
$$I^+(n, (2k-1,4k))= 
I(\underbrace{1, \dots ,1}_{n-k-2},2,\underbrace{1, \dots ,1}_{k},2,
\underbrace{1, \dots ,1}_{k},2,\underbrace{1, \dots ,1}_{n-k-2}).$$

\item Letting $s_1= 2,2,1,2,1,2,1,2$ and $s_2 = 2,1,2,1,2,1,2,1$, 
 then $$I^+(n,(1,4))= I^-(n,(1,3))=$$
$$\begin{cases}
I(\underbrace{s_1,s_2, \dots}_m ) & \text{if $n= 5m$,}\\
I(1,2,\underbrace {s_2, s_1, \dots}_m )& \text {if $n= 5m+1$,}\\
I(1,1,2,1,\underbrace {s_1,s_2, \dots}_m) & \text {if $n = 5m+2$,}\\
I(2,1,2,1,2,\underbrace {s_2,s_1, \dots}_m) & \text {if $n= 5m+3$, and}\\
I(1,2,1,2,1,2,1,\underbrace {s_1,s_2, \dots}_m) & \text {if $n = 5m+4$}
\end{cases} 
$$

\item Letting $s_1= 2,1,2,1,2$, $s_2= 1,2,1,2,1$, 
$s_3= 2,1,2,1,1$, $s_4 = 1,2,1,1,2$ and $s_5 = 2,1,1,2,1$, then 
$$ I^+(n,(1,3))= I^-(n,(3,8))=$$
$$\begin{cases}
I(\underbrace{s_1, s_2, s_3 , s_4, s_5,  \dots}_m ) & \text{if $n= 3m$,}\\
I(1,2,\underbrace {s_2, s_5, s_3, s_1, s_4 \dots}_m )
& \text {if $n= 3m+1$, and}\\
I(1,1,2,1,\underbrace {s_3, s_1, s_4, s_2, s_5, \dots}_m) 
& \text {if $n = 3m+2$}
\end{cases}$$

\end{enumerate}

The first item holds for $n=1$ by Proposition~\ref{endpt}. 
By Proposition~\ref{mult} and induction on $n$,
it holds for all $n$.  

By Proposition~\ref{slopes}, $I^-(n,(1,2))$ contains 
$I^-(n-1,(1,2))I^-(1,(1,2))$.  By induction on $n$, this is the ideal 
$$I(\underbrace{1, \dots ,1}_n,2,\underbrace{1, \dots ,1}_{n-2},2).$$
Proposition~\ref{slopes} and (1) imply that $I^-(n,(1,2))$ contains
$y^{2n+1}$. Checking colengths, this forces $$I^-(n,(1,2))=
I(\underbrace{1, \dots ,1}_n,2,\underbrace{1, \dots ,1}_{n-1}).$$

The third item follows from Proposition~\ref{mult} after checking the first
three ideals by hand via Proposition~\ref{endpt} and checking that all
of the ideals have the correct colength.  This can be done by
induction, checking the cases of the 
three possible $r$'s separately.  

For $n\le 3$, (4) follows from 
Proposition~\ref{endpt}.  Then by induction on $n$ and checking
colengths, by Proposition~\ref{mult} we have 
$$I^{+}(n,(0,1))=I^{+}(n-2,(0,1))I^{+}(2,(0,1))
+I^{+}(n-3,(0,1))I^{+}(3,(0,1))$$ 
which is the claimed ideal.

By Proposition~\ref{endpt}, in all characteristics except possibly
$2$, $$I^+(1, (0,-1))= I^-(1, (1,2))= I(1,2).$$  Since the polytope
in the construction of a standard fan for positive
characteristic must be contained in the corresponding polytope in
characteristic $0$, this last result holds also in characteristic $2$ and
the boundary of $F(I(4))$ is as described in Figure~\ref{boundary} in 
all characteristics.   

We will prove (5) and (6) by induction on $n$. For $n \ge 2$ by 
Proposition~\ref{mult} we have 
(5) for  $\left [ \frac {n+1}{2} \right ] \le k \le n-2$ and (6) for 
$\left [ \frac {n-1}{2} \right ] \le k \le n-3$ 
since for 
$$ \left [\frac {n+1}{2} \right ] \le k <n-2 $$ 
we have 
$$I^{-}(n,(2k-1,4k))=
I^{+}(n,(2k-3,4k-4))= $$
$$I^{-}(n-1,(2k-1,4k))I^{-}(1,(2k-1,4k))+$$
$$I^{-}({\left[\frac {n}{2}\right]},(2k-1,4k))
I^{-}({\left[\frac {n+1}{2}\right]},(2k-1,4k))$$
$$=I(\underbrace{1,\dots ,1}_{n-k-1},2,\underbrace {1,\dots, 1}_{k-1},2,
\underbrace{1,\dots ,1}_{k-1},2,\underbrace{1,\dots ,1}_{n-k-1}).$$ 

Thus for $n\ge 3$, 
$$I^{-}(n,(2n-5, 4n-8))=  
I(1,2,\underbrace {1,\dots, 1}_{n-3},2,
\underbrace{1,\dots ,1}_{n-3},2,1).$$ 
By Proposition~\ref{mult}, $I^{+}(n,(2n-5, 4n-8))$ is contained in 
$I^+(1,(2n-5, 4n-8))I^{+}(n-1,(2n-5, 4n-8))$.
Proposition~\ref{slopes} then forces $I^{+}(n,(2n-5, 4n-8))$ to be as
claimed.  By Proposition~\ref{mult} the points corresponding to 
$2$-dimensional cones in the standard fan $\Delta ((x,y^4)^n)$ lying
in the convex cone bounded by the rays through $(1,2)$ and
$(2n-5,4n-8)$ correspond to ideals contained in 
$$J=I^-(n-1, (1,2))I^-(1,(1,2))= 
I(\underbrace{1,\dots ,1}_{n},2,\underbrace {1,\dots, 1}_{n-2},2).$$
The cones corresponding to the ideals $I^-(n, (1,2))$ and 
$I^+(n,(2n-5, 4n-8))$ lie in the cone bounded by the rays through $(1,2)$ and
$(2n-5,4n-8)$.  Suppose there is another such cone corresponding to an
ideal generated by $x^ey^f$ over $J$.
Since $I^-(n, (1,2))$ and 
$I^+(n,(2n-5, 4n-8))$ are generated over $J$ by $y^{2n+1}$ and 
$x^{2n-1}$ the inequalities 
$$ \frac {2(2n-1-e) -f}{3(2n-1-e) -f} > 
\frac{2e - 2n-1+f}{3e - 2n-1+f} > \frac{1}{2}$$
follow from Corollary~\ref{slopecor}.  Thus it follows that $e+f = 2n$
and $e \ge f$.  However, this implies that $x^ey^f$ is already
contained in $J$.  Therefore we have arrived at a contradiction and 
by Corollary~\ref{slopecor}, 
$I^+(n,(2n-5, 4n-8))= I^-(n, 2n-3, 4n-4)$ and 
$I^-(n, (1,2))= I^+(n, 2n-3, 4n-4)$.  If $n= 2m+1$, the ideal 
$I^-(n, (2m-1, 4m))$ cannot be deduced from Proposition~\ref{mult} 
alone, containment in $I^-(n-1,(2m-1, 4m))I^-(1, (2m-1, 4m))$ together
with knowledge of $I^+(n, (2m-1, 4m))$ and Proposition~\ref{slopes}
force it to be as claimed.  

From what we have shown so far, it follows that the boundary of
$\bar F((x,y^4)^n)$ for $n\le 4$ is as depicted in Figure~\ref{boundary}.  
Each vertical line segment 
represents the boundary divisor corresponding to the $1$-dimensional 
cone passing through the lattice pointed as marked.  The ideals in 
between the line segments represent the intersections of the 
corresponding boundary divisors.

We will now find the boundary of $\bar F((x,y^4)^5)$.  By 
Proposition~\ref{mult} $I^+(5,(1,4))$ is contained in
$I(1,1,2,1,2,1,2,1,2)$.  Thus Proposition~\ref{slopes} forces 
$I^+(5,(1,4))= I(2,2,1,2,1,2,1,2)$.  In the same manner that 
Corollary~\ref{slopecor} was used to deduce $I^-(n,(1,2))= 
I^+(n, (2n-3, 4n-4))$, it can be used to deduce that 
$I^+(5,(1,4))= I^-(5,(1,3))$ and $I^-(5,(3,8))= I^-(5,(1,3))$, 
completing the diagram for $\bar F((x,y^4)^5)$.  

It remains to verify the last two items.  These items follow  
from induction on $n$ and Proposition~\ref{mult} since 
$$I^+(n,(1,4))=
I^+(n-2,(1,4))I^+(2,(1,4))+I^+(n-5,(1,4))I^+(5,(1,4))$$
and 
$$I^+(n,(1,3))=
I^+(n-3,(1,3))I^+(3,(1,3))+I^+(n-5,(1,3))I^+(5,(1,3)).$$

Thus we have demonstrated our eight claims, proving the
boundaries of the spaces $\bar F((x,y^4)^n)$ to be as in
Figure~\ref{boundary} for $n \le 6$.  
We stop here, not because our techniques do
not suffice for larger $n$, but because we have illustrated them all.
  
\section{Almost Toric Varieties}

In the last section we developed an arsenal of techniques for 
understanding the spaces $F(I)$ for ideals $I$ with certain 
measuring sequences.  The 
basic idea was to use the structure of the normalization of $F(I)$ as
a toric variety together with the correspondence of the points of
$F(I)$ as ideals.  In this section we take a baby step towards
developing a similar arsenal for dealing with other measuring
sequences by considering a few examples.  In these examples 
we consider some spaces 
$\bar F(I_1, \dots ,I_r)$ in which the $I_j$'s have
measuring sequence at most $m(5,1)$.  These spaces will be three
dimensional, with a two dimensional torus acting on them.  In this
sense, they are almost toric varieties.  This observation will 
enable us to associate data to the boundary divisors analogous to 
the rays in the standard fan of the previous section.  

For the remainder of this section we will use $I_j$ to denote the
ideal $(x,y^j)$.

\begin{example}\label{J5}
\end{example}
Consider the space $F(I_3,I_4,I_5)$.  There is a natural
projection from this space to the space $F(I_3,I_4,I_5,I_1I_4)$ 
which is the projectivization of the bundle $V(I_4/I_1I_4)$ (as 
defined in \cite {R3}) over 
$F(I_3,I_4, I_1I_4)$.  Roughly $V(I_4/I_1I_4)$ is the bundle 
with fiber $J_2/J_3$ over a point 
$(J_1, J_2, J_3)\in  F(I_3,I_4,I_1I_4)$.  
This projection contracts the $\P^1$ of points of the 
form $$(I_1^2,I_1I_2,I_2^2+I_1^3,(sx^2+ty^3)+I_1^2I_2)$$ 
to the point 
$$c = (I_1^2,I_1I_2,I_2^2+I_1^3)$$ 
and is a local isomorphism everywhere else.  Although the $\P^1$ lying 
above $c$ separated the divisors $D(0,1)$ and 
$D(1,2)$, corresponding to the rays through $(0,1)$ and $(1,2)$
respectively in the relevant standard fans, the point 
$c$ lies on both of these divisors.  From Table 2 \cite{R3}, 
we see that the preimages of $D(0,1)$ and $D(1,2)$ in the fiber of 
$J_5$ over $P$ 
are the loci of points of the form 
$$((sx+ty^2,y^3),I_1I_2, (uxy+vy^3,x^2,xy^2,y^4))$$
and $$(I_1^2, (sxy+ty^2,x^2,xy^2,y^3),(u(sxy+ty^2)+vx^2)+I_1^3)$$ 
respectively.  
Since these two boundary divisors lie in planes in the 
Pl{\"u}cker embedding meeting only at $c$, 
the Zariski tangent space to $c$ in $F(I_3,I_4,I_5)$ 
must have dimension at least $4$.
Hence the singular locus of $F(I_3,I_4,I_5)$ consists of the
point $c$.     

\vspace {.1 in}
 
\begin{example}\end{example}

\begin{figure}\label{P} \hspace{3 in}
\begin{picture}(-100,0)(300,200)
\hspace{3 in}

\put(0,0){\circle*{9}}
\put(20,0){\circle*{9}}
\put(40,0){\circle*{9}}
\put(60,0){\circle*{9}}
\put(80,0){\circle*{9}}
\put(100,0){\circle*{9}}
\put(120,0){\circle*{9}}
\put(140,0){\circle*{9}}
\put(160,0){\circle*{9}}
\put(180,0){\circle*{9}}
\put(200,0){\circle*{6}}
\put(220,0){\circle*{6}}
\put(240,0){\circle*{6}}
\put(260,0){\circle*{3}}
\put(280,0){\circle*{3}}
\put(300,0){\circle*{3}}

\put(0,20){\circle*{9}}
\put(20,20){\circle*{9}}
\put(40,20){\circle*{9}}
\put(60,20){\circle*{9}}
\put(80,20){\circle*{9}}
\put(100,20){\circle*{9}}
\put(120,20){\circle*{9}}
\put(140,20){\circle {9}}
\put(160,20){\circle*{6}}
\put(180,20){\circle*{6}}
\put(200,20){\circle*{6}}
\put(220,20){\circle*{3}}
\put(240,20){\circle*{3}}
\put(260,20){\circle*{3}}

\put(0,40){\circle*{9}}
\put(20,40){\circle*{9}}
\put(40,40){\circle*{9}}
\put(60,40){\circle*{9}}
\put(80,40){\circle {9}}
\put(100,40){\circle*{6}}
\put(120,40){\circle*{6}}
\put(140,40){\circle*{6}}
\put(160,40){\circle*{6}}
\put(180,40){\circle*{3}}
\put(200,40){\circle*{3}}
\put(220,40){\circle*{3}}

\put(0,60){\circle*{6}}
\put(20,60){\circle*{6}}
\put(40,60){\circle*{6}}
\put(60,60){\circle*{6}}
\put(80,60){\circle*{6}}
\put(100,60){\circle*{6}}
\put(120,60){\circle {6}}
\put(140,60){\circle*{3}}
\put(160,60){\circle*{3}}
\put(180,60){\circle*{3}}

\put(0,80){\circle*{6}}
\put(20,80){\circle*{6}}
\put(40,80){\circle*{6}}
\put(60,80){\circle {6}}
\put(80,80){\circle*{3}}
\put(100,80){\circle*{3}}
\put(120,80){\circle*{3}}
\put(140,80){\circle*{3}}

\put(0,100){\circle*{3}}
\put(20,100){\circle*{3}}
\put(40,100){\circle*{3}}
\put(60,100){\circle*{3}}
\put(80,100){\circle*{3}}
\put(100,100){\circle {3}}

\put(0,120){\circle*{3}}
\put(20,120){\circle*{3}}
\put(40,120){\circle {3}}

\end{picture}

\vspace{3 in}
\caption{$F(3;4;1,4;5;1,5)$}
\end{figure}
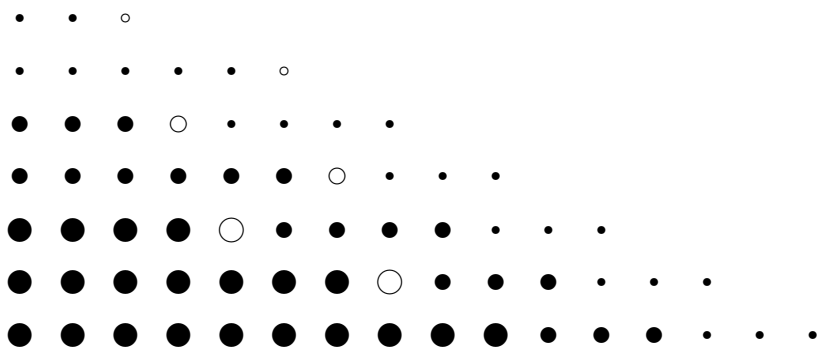 

Consider the space 
$$Y= F(I_3,I_4, I_1I_4, I_5, I_1I_5).$$  By Lemma 3.3 of [10], 
the coset representatives
for $G(I_2)$ under right multiplication by 
$G(I_2,I_3,I_4, I_1I_4, I_5, I_1I_5)$ are automorphisms $g$ with 
$$g(x)= x+ay^2+by^3+cy^4$$
and 
$$g(y)=y.$$  Thus $a$, $b$, and $c$ form coordinates on the interior
of $Y$. 
The picture 
$P$ obtained by plotting the points of the exponent vectors of the 
non-zero terms of the coordinate functions for $Y$ is shown in
Figure~\ref{P}.  The $x$ and $y$ coordinates of the dots
correspond to the exponents of $a$ and $b$ respectively.  
The sizes of the dots correspond to the exponent of $c$.  Open circles
are used for points for which the corresponding monomial is not in the
span of the coordinate functions for $Y$.  The generators for the
vector space spanned by the coordinate functions modulo the space
spanned by the monomials in the span are 
$$a^5(ac-b^2)(ac-2b^2)$$ and 
$$a^2b(ac-b^2)(ac-2b^2).$$  
Henceforth, consider $Y$ to be embedded in projective space with 
a basis of the span of the coordinate functions consisting of 
monomials together with these two functions as coordinate functions.  
There are $8$ faces of the convex hull $H$ of $P$, five corresponding
to subvarieties of the boundary and three corresponding to the 
closures of the coordinate planes.  Of the latter five planes, all but
the one passing through the open circles correspond to boundary
divisors.  The one passing through the open circles corresponds to the
line of points with all coordinates with monomial coordinate functions
zero.  There are two more boundary divisors, corresponding to the two 
roots of the non-monomial coordinate functions.  
Thus, the point whose only nonzero 
coordinate is the one with coordinate function $a^5(ac-b^2)(ac-2b^2)$
is singular.

\end{document}